\documentclass[11pt,english]{article}
\usepackage[dvips]{graphics}
\usepackage{graphicx}
\usepackage{verbatim}
\usepackage[english]{babel}
\usepackage{multicol}
\usepackage{color}
\usepackage{amsmath,amssymb,dsfont}
\usepackage[utf8]{inputenc}

\usepackage{hyperref}
\usepackage{fancyhdr}
 \font\ninerm=cmr9
\long\outer\def\abstract#1{\bigskip\vbox{\noindent\ninerm
\baselineskip=10pt#1}\nobreak\bigskip}

\newcommand{\R}{\mathbb{R}}

\newcommand{\E}{\mathbb{E}} 

\newcommand{\V}{\mathbb{V}}
\newcommand{\indic}[1]{{\mathbb I}_{#1}}
\newcommand{\proba}{\mathbb{Q}}

\newcommand{\Fr}{\mathcal{F}}

\newcommand{\Cr}{\mathcal{C}}
\newcommand{\Ir}{\mathcal{I}}
\newcommand{\Ur}{\mathcal{U}}
\newcommand{\Ar}{\mathcal{A}}

\newtheorem{definition}{Definition}
\newtheorem{proposition}{Proposition}
\newtheorem{lemma}{Lemma}
\newtheorem{rem}{Remark}
\newtheorem{corollary}{Corollary}

\newtheorem{thm}{Theorem}

\begin{document}
\title{The Design of Optimal Re-Insurance Contracts when Losses are Clustered}

\author{Guillaume Bernis\thanks{BPCE Vie. The views expressed in this paper reflect those of the authors and not necessarily those of BPCE Vie, email:  \texttt{guillaume.bernis@bpce.fr}.}\and Cristina Di Girolami\thanks{corresponding author, Dipartimento di Matematica, Universit\`a Alma Mater Studiorum Bologna, Piazza di Porta San Donato, 5, 40126 Bologna BO, Italy, email: \texttt{cristina.digirolami2@unibo.it}.} \and Simone Scotti \thanks{Dipartimento di Economia e Management, Universit\`a di Pisa, via Ridolfi 10; email:   \texttt{simone.scotti@unipi.it }.}}
\setlength{\parskip}{0mm}
\setlength{\parindent}{0mm} \setlength{\hoffset}{0cm}
\setlength{\textwidth}{14cm}

\maketitle

{\footnotesize
{\bf Abstract} {This paper investigates the form of optimal reinsurance contracts in the case of clusters of losses. The underlying insured risk is represented by a marked Hawkes process, where the intensity of the jumps depends not only on the occurrence of previous jumps but also on the size of the jumps, which represents the financial magnitude of the loss. The reinsurance contracts are applied to each loss at the time of occurrence, but their structure is assumed to be constant. We derive closed-form formulas within the mean-variance framework. Additionally, we demonstrate that the optimal contract is not the classical excess-loss (deductible) form. The optimal contract is piecewise linear with three ranges: first, no reinsurance below a certain threshold; second, reinsurance with a slope greater than $1$; and finally, full reinsurance. When the marked process converges to a Poisson process, we recover the optimality of the deductible form.}
}
\medskip

{\bf Keywords}: Optimal Reinsurance; Optimal Contract; Hawkes processes; Clusters.

{\bf JEL Classification}: G22, C61, G32.

{\bf AMS Classification}: 91G05, 60G55.

\newpage 

%
%
%
%
%
%

\section{Introduction}

The problem of optimally setting up an insurance policy is an old and pivo-tal question for insurance and reinsurance companies, as highlighted in the seminal paper by Borch \cite{Bor}. A key result, due to Arrow \cite{Arrow}, is that the optimal contract for a risk-averse buyer {\it takes the form of 100 percent coverage above a deductible minimum}; this is known as the {\it excess-of-loss contract}. Arrow's analysis is quite general (see possible extensions in \cite{GolSch} and references therein): he assumes that the premium is a proportion of the actuarial value of the policy and that the insured chooses the coverage function based on standard utility maximization. The payments are non-negative, and the description is static, focusing on a one-period model with a premium and a random loss.
The extension to a continuous-time framework was achieved by Aase \cite{Aas}, where the occurrence of losses is represented by a compound Poisson process, and the insurance premium is paid continuously over time. Arrow's result is preserved essentially because of the infinite divisibility with respect to time. According to \cite{Aas} and as extended by \cite{EecGolSch} and \cite{GolSch}, the formulation of the criterion remains purely static, applying the reinsurance contract to a single aggregated loss over a fixed time window. The model thus belongs to the Cram\'er-Lundberg setup; see also Bauerle \cite{Bauerle} and Touzi \cite{Touzi}.
However, while the premium payment can be reasonably assumed to be proportional to time, since it is paid on a fixed schedule, the hypothesis of a constant accident arrival rate per unit of time is too restrictive, as noted by Grandell \cite[Chapter 2]{Grandell}. Losses are observed, or can be considered, occurring in clusters in many insurance areas, such as lapse risk in life insurance \cite{Barsotti},  cyber risk \cite{BalGhe, BouCherHil1}, natural disaster for insurance  
\cite{LesDeaLej} and mortality intensity \cite{JiZh}. 

In recent years, attempts have been made to relax the Poisson framework by using Cox processes; see, for instance, Albrecher and Asmussen \cite{AlAs06}, Brachetta and Ceci \cite{BraCec}, Dassios and Jang \cite{DassJan} and Embrechts et al. \cite{EmSchGra}.
The main issue with Cox processes is that there is a hidden process describing the intensity of accident arrivals, which needs to be inferred from data. There is a substantial and growing body of literature on filtering in insurance; see, for instance, \cite{BraCec2, DassJan2}.
In order to overcome this difficulty and to maintain a parsimonious model, we model the arrival of losses using a Hawkes process, as discussed in \cite{Hawkes}. This resulting model is fully observable since the times of loss arrivals are known as in all the standard models in insurance. In our set-up, the intensity will increase at each of these times and will exhibit exponential decay between accidents. The use of a Hawkes process to model accident arrivals has garnered interest from researchers in recent years, including Dassios and Zhao \cite{DassZh}, Garetto et al. \cite{GarLeoTor}, but also in cyber risk \cite{BouCherHil1, HilRevRos, CallFonHillOng}, natural disasters \cite{LesDeaLej}. Hawkes processes can be viewed as a parsimonious form of a Cox process, see Albrecher and Asmussen \cite{AlAs06}, in which the intensity is driven by the jumps themselves.

However, to our knowledge, the shape of the optimal contract has not yet been studied in this framework, as previous works have directly considered the cases of proportional and/or excess-of-loss reinsurance.
We point out that, in the presence of accident clusters, the arguments in \cite{Aas, EecGolSch, GolSch} to reduce the problem to the static one {\it \`a la Arrow} are no longer valid. When focusing on a cluster of accidents, the insurance buyer will not be covered by the seller at the deductible threshold for each loss. The insurance buyer is, therefore, concerned about the increasing frequency of losses and may prefer a contract with a different shape.

To study this problem, we propose a marked Hawkes model for loss arrivals, as discussed in \cite{BerSalSco18}. The marked Hawkes setup captures both the clustering effects and the self-exciting feature: a significant accident is more likely to trigger the arrival of new accidents compared to smaller ones. The marked Hawkes process, along with its generalization to branching processes, has been extensively utilized in recent financial literature to capture the evolution of stochastic volatility and the VIX index \cite{HorstXu2, JMSZ}, commodities \cite{BriGonSga, CalMazSga}, term structures \cite{BGSS23, JMS}, and market microstructure \cite{AberJedi, BacMuz, HorstXu}. The case of non-exponential decay, which represents the non-Markovian scenario, has been more recently studied using stochastic Volterra equations with jumps; see \cite{AbiLarPul, BLP24, BPS22}.

To the best of our knowledge, our work is the first to focus on the problem of optimal reinsurance contract design within a marked Hawkes framework, where accidents can be clustered and the magnitude of a loss has a direct impact on the occurrence of future accidents. We will study the case of a generic feedback effect and then specify it into constant and linear impact, that is the two more standard cases. The first where the excess of the intensity after a loss does not depend on the loss size itself. The second, that is more natural, 
where the large losses have more impact on the future losses that the small ones.

\begin{figure}[htbp]
	\centering
		\includegraphics[width=0.50\textwidth]{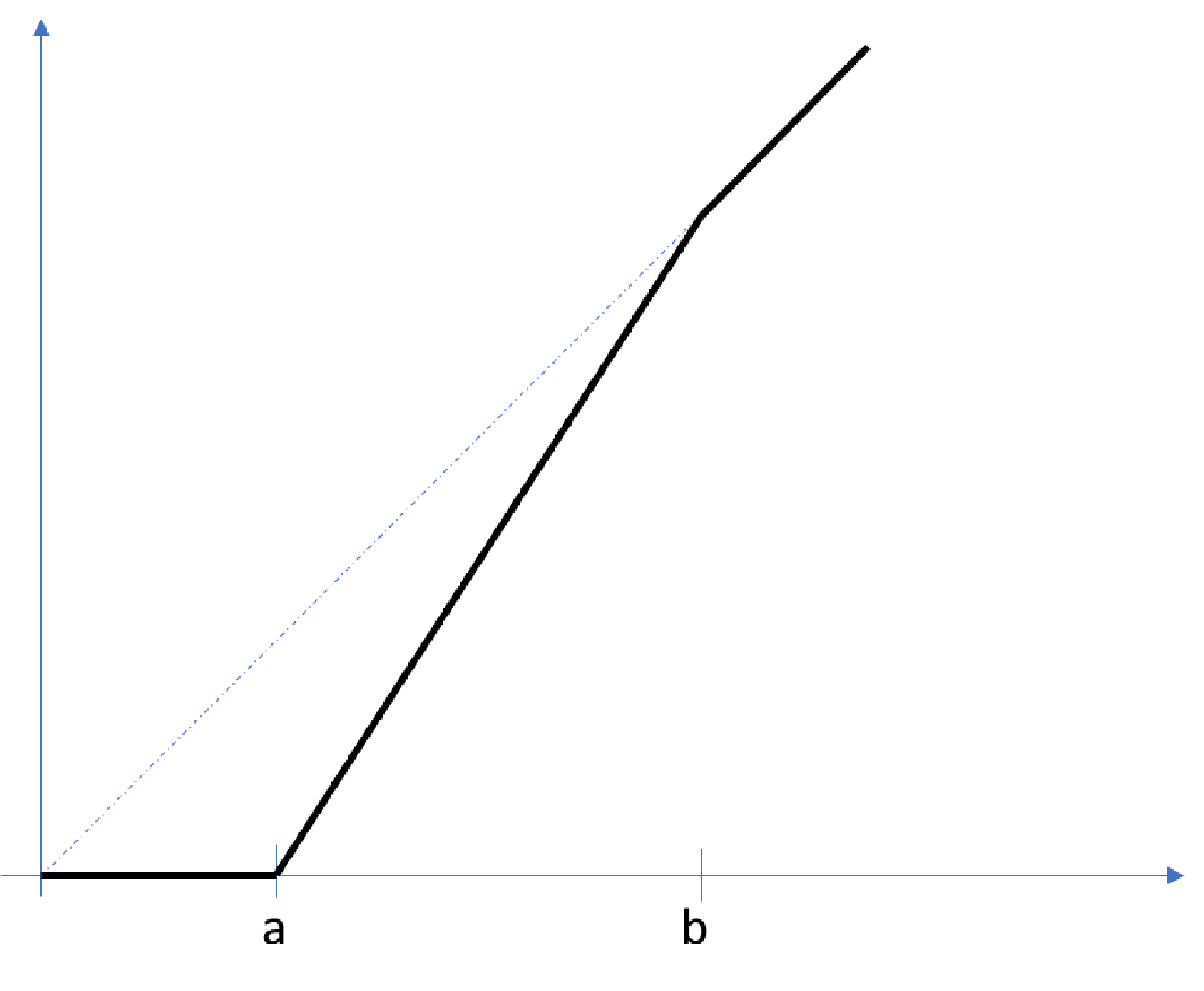}
\begin{minipage}{0.5\textwidth}
	\caption{Optimal reinsurance contract. Claim size on $x$-axis and covered amount on $y$-axis.}
	\label{fig:OptimalSolution}
\end{minipage}
\end{figure}

Our main result is that the optimal contract is no longer the deductible one in the second setup. Instead, it is piecewise linear with three ranges. First, as in the deductible case, there exists a threshold $a$ below which there is no reinsurance. Conversely, there is a second threshold
 $b>a$  above which the optimal contract provides total reinsurance. Finally, between the two thresholds $[a,b]$, the reinsurance is affine with respect to the loss, ensuring the continuity of the reinsurance contract relative to the loss; see Figure \ref{fig:OptimalSolution}, where, on the $x$-axis, there is the claim size and, on $y$-axis, the covered amount. This result is obtained under the assumption of a mean-variance utility function, in accordance with Bauerle \cite{Bauerle}. 

At first glance, one might consider that because the slope of the optimal contract is greater than $1$, this strategy could induce moral hazard, in line with the analysis of Huberman et al. \cite{HubMaySmi}. However, ex-ante moral hazard (specifically, the distortion of the insurer's dynamic risk) does not arise, as the reinsurance payment is lower than the loss incurred. Furthermore, the intensity of loss occurrence increases after each event, which means that the insurer remains risk-averse even in case of very large losses (that are perfectly covered by the contract), since such losses could trigger a cascade of subsequent events. 

For instance, in the case of health insurance, classical theory postulates that having insurance coverage diminishes an individual's incentive to engage in preventive measures to maintain their health, see Harris and Raviv \cite{HarRav}. However, this theory is based on the implicit assumption that the occurrence of new diseases is independent of prior incidents, a strong hypothesis that can be easily falsified. An example of this can be found in the application of Hawkes processes in health management \cite{GarLeoTor}.

Focusing on ex-post moral hazard (specifically, the distortion of claims), this issue also does not arise in our model. In fact, an increase in the claim size leads to a more significant shift in the intensity of new events. This outcome is a result of the marked Hawkes setup we have chosen, in contrast to a standard Hawkes process, where the intensity is influenced solely by the occurrence of events.
For example, in the context of cyber insurance, an insurer might delay declaring a vulnerability, thereby increasing the associa-ted claim. However, postponing the official declaration, also provides more opportunities for cyber attacks, delays the updating of information systems, and ultimately heightens the risk of a cyber pandemic. 

The rationale behind our result is that the reinsurance buyer pays a constant premium per unit of time but is potentially exposed to the presence of clusters. Moreover, this risk is magnified by the impact of the marks on the intensity of future losses. The reinsurance buyer is risk-averse and prefers to be fully covered for significant accidents, as they will need to account for potential future losses.
At first glance, the particular shape of the optimal contract appears to transfer completely the risk of large accidents to the reinsurer. However, this analysis neglects the feedback effect resulting from the cascade of small losses induced by the initial large accident. This effect has been highlighted in the Hawkes literature as a cluster decomposition; see Hawkes and Oakes \cite{HawOak}, as well as more recent works \cite{HorstXu2, JMSZ}. 

After this Introduction, in Section \ref{secDynReinFram}, we present the stochastic representation of the insured risk through a marked Hawkes process. We define the insurer's wealth process and the reinsurance contracts. Finally, we introduce the mean-variance criterion to be considered in the decision-making process and provide an explicit formula. In Section \ref{secOptReinCont}, we prove the existence of an optimal reinsurance contract and outline its form. We also investigate the limiting case when the Hawkes process converges to a Poisson process and demonstrate that we recover the deductible form.

\section{The Dynamic Reinsurance Framework}		\label{secDynReinFram}

\subsection{Stochastic representation of the risk}			\label{secStoRepRis}

Let us set-out the following stochastic framework. Let $(\Omega, \{\Fr_t\}_{t\geq 0}, \Fr, \proba)$ be a filtered probability space, satisfying the usual conditions, see for instance
 \cite[page 10]{KarShr91},
equipped with a  marked Hawkes process, represented by its Poisson measure $\nu(dt,dz)$. Hereafter, we define the marked Hawkes process following \cite[Definition 2]{BerSalSco18} and \cite{BriSga20}.

\begin{definition}\label{defHawkes}
The point process $\{\lambda_t\}_{t\ge 0}$  is a marked Hawkes process if its compensator is given by $\Theta(dz) \lambda_{t^-} dt$, where $\Theta(dz)$ is a measure on $(\R^+, {\mathcal B}(\R^+))$, and $\lambda_t$ reads
\begin{equation}\label{lambda-SDE}
{\lambda}_t  = \lambda_0 + \beta \int_0^t(\underline{\lambda} -  \lambda_{s^-})ds + \int_0^{t^-}\int_{\R^+} f(z)\nu(ds,dz)
\end{equation}
with $\lambda_0 > \underline{\lambda}$,  the long run intensity  and $f$ is a given non-negative function, square-integrable with respect to the measure $\Theta(dz)$. 
Moreover we assume the following integrability assumption: 
\begin{equation} 		\label{eqINTtheta}
\int_0^{+\infty}z^2\Theta(dz)<+\infty.
\end{equation} 
\end{definition}

This function $f$ captures the self-exciting effect. We will be interested in two polar cases:\begin{description}
\item[Constant impact:] the function $f$ is a constant equal to $\overline{f}$. This is the case where the excess of intensity does not depend on the size of the previous claims
\item[Linear impact:] $f(z) = \Lambda z$, for some non-negative constant $\Lambda$. This is the central case of our main result. From a practical point of view, it is more natural to assume that the large losses have more impact on the future losses that the small ones.
\end{description}

Throughout the paper, we assume that $\Theta(dz)$ admits at least moments up to order two. 
We recall that $\nu(ds,dz)$ is the Poisson measure associated to the $(T_i,Z_i)_{i\ge 1}$, where the $T_i$ are the times of occurrence and the $Z_i$ the sizes of the jumps. In the following, each $T_i$ will represent the time of occurrence of an insured loss and $Z_i$ will be the financial cost associated to this loss.
From now on, we assume the following usual mean-reverting condition for norm of the kernel, which can be reformulate as the $L^1(\Theta)$-norm smaller than the mean-reversion parameter $\beta$, i.e.
\begin{equation}
\displaystyle   \int_{\R^+}f(z)\Theta(dz)< \beta. \label{hypErgodic} 
\end{equation}
Under assumption \eqref{hypErgodic}, the Hawkes process $\{\lambda_t\}_{t\ge 0}$ is well defined, mean-reverting, admits moments up to order two and its Laplace Transform is known, see \cite[Proposition 7.3, p. 176]{BerSco20}. For sake of readability, set out the following functions of  time: 
\begin{equation} 
 m(t) := \E\left[ \lambda_t\right] , \quad  M(t) := \int_0^t m(s) ds = \E\left[\int_0^t\lambda_s ds\right],
\end{equation} 
\begin{equation} 
 m^{(2)}(t) := \E\left[\lambda^2_t\right],  \quad \mbox{ and }\quad  M^{(2)}(t) := \int_0^tm^{(2)}(s)ds=\E\left[\int_0^t\lambda^2_s ds\right].
\end{equation} 

\begin{definition}\label{defH} Let $\mathcal{H}$ denote the following set of functions 
$$
\mathcal{H}:= \left\{ h: \R^+ \longrightarrow \R \ ; \mbox{ such that $\Theta$-a.s. } \vert h(z) \vert \le z \right\}.
$$
For every $h\in \mathcal{H}$ 
we define the process $\{X_t[h]\}_{t\ge 0}$ by
\begin{equation}\label{eqXh}
X_t[h] := \int_0^t\int_0^{+\infty}h(z)\nu(ds,dz) 
\end{equation}
and the functional $H: \mathcal{H}\longrightarrow \R$ defined as
\begin{equation}	\label{defHfunctional}
\displaystyle H[h]:= \int_0^{+\infty}h(z)\Theta(dz).
\end{equation}
In the sequel we will denote by $\Ir$ denotes the identity function on $\mathbb{R}^+$.
\end{definition}

Next proposition gives the expected value of the process $X[h]$.

\begin{proposition}\label{Prop1}
Let $h\in \mathcal{H}$, then the process $X[h]$ is integrable and its expectation reads
\begin{equation}\label{eqEspX}
\E\Big[X_T[h] \Big] = 
  \int_0^{+\infty}h(z)\Theta(dz) M(T) = H[h] \cdot M(T)
\end{equation}
\end{proposition}

{\bf Proof:}
We first remark that $\vert X_T[h] \vert \leq X_T[\Ir]$.
As a consequence, the integrability of $X_T[\Ir]$ guarantees the integrability for all $X_T[h]$. 
It is now easy to note that $X_T[\Ir]$ is the cumulated loss process: see the last term in equation 
\eqref{lambda-SDE}. The integrability of the cumulated loss process is guaranteed by the 
square integrability assumption on $\Theta$ given in \eqref{eqINTtheta}, see for instance \cite[Section 5]{BerSco20}.
We now focus on the computation of the expectation, we have
$$
\E\Big[X_T[h] \Big] = \E\left[ \int_0^t\int_0^{+\infty}h(z) \left\{ \nu(ds,dz) - \lambda_{s-} \Theta(dz) ds +  \lambda_{s-} \Theta(dz) ds \right\} \right]
$$
Since $h\in \mathcal{H}$ and the integrability condition obtained before, we have that 
$\int_0^t\int_0^{+\infty}h(z) \left\{ \nu(ds,dz) - \lambda_{s-} \Theta(dz) ds  \right\}$ is a martingale. The final result is then obtained by Fubini theorem. 
\hfill $\Box$

\medskip
Another possible approach is to exploit the link between Poisson and Hawkes process, see
\cite[Proposition 2]{BerSalSco18} for the case of a bounded jump size, 
 \cite[Section 2]{BraCalCecSga} for the unbounded with technical integrability condition and the condition can be relaxed to include our case following \cite{Yan}.

\subsection{The insurer program}

In this section, we express the wealth of an insurer that can access to a reinsurance contract in order to mitigate the global risk exposure. Our formulation follows the setup {\it \`a la Cramer-Lundberg} developed in Bauerle \cite{Bauerle}. See Brachetta et al. \cite{BraCalCecSga} for an adapted setup with Hawkes driver. In particular, we consider that the reinsurance contract applies at any loss $Z_i$ occurring at any time $T_i$, not on an aggregated wealth at a final time.

The insurer is endowed with a certain capital amount $R_0$ at time $t=0$. It is assumed to insure all the claims generated by the marked Hawkes process. 

The premium of these claims is assumed to be paid continuously over time with a rate $\rho   \overline{\theta}$, expressed in Euros by year. Quantity $ \overline{\theta}$ denotes
\begin{equation}
 \overline{\theta}:=\displaystyle   \int_{\R^+}z\Theta(dz)=H[\Ir]
\end{equation}
and represents the average loss associated to the sinister. We observe that, by Definition \ref{defHfunctional} the mean-reverting assumption \eqref{hypErgodic} rewrites as $H[f]<\beta$.

\begin{definition}\label{defReinsContract} A reinsurance contract is a function in ${\mathcal H}$, non-negative $\Theta$-almost everywhere. We denote by $\Cr$ the set of reinsurance contracts, i.e. $\Cr\subset \mathcal{H}$ and 
$$
\Cr:=\left\{ h: \R^+ \longrightarrow \R \ ; \mbox{ such that $\Theta$-a.s. }  0\le  h(z)  \le z \right\}.
$$
\end{definition}

The cost of such a contract is assumed to be linear in $\phi$, and paid continuously over time with a rate
$$ c\int_{0}^{+\infty}\phi(z)\Theta(dz).$$

The total wealth of the insurer at time $t$, applying reinsurance contract $\phi$, is given by
\begin{equation}\label{eqRinitial}
R_t(\phi) = R_0 + \int_0^t\int_{0}^{+\infty}\left( \phi(z)-z \right) \nu(ds,dz) + \rho\overline{\theta}t-c t\int_{0}^{+\infty}\phi(z)\Theta(dz) \, .
\end{equation}

With the notations given in Definition \ref{defH}, we can rewrite the insurer's wealth at final horizon time $T$ of the program as
\begin{equation}\label{eqRTphi}
R_T(\phi) = R_0 + (\rho-c)T \cdot \overline{\theta} + X_T[\phi - \Ir] - c T H[\phi-\Ir] .
\end{equation}
We observe that the stochastic part is located only in the third term of the right hand side.

We consider the following mean-variance criterion, in the spirit of Bauerle \cite{Bauerle}:
\begin{equation}	\label{defMVcr}
\Ur(\phi) := \E\left[R_T(\phi)\right] - \gamma \V\left[R_T(\phi)\right]\, ,
\end{equation}
where $\gamma$ is a positive constant and $\V$ denotes the variance. 


%

\begin{thm}\label{thmMeanVar} The mean-variance criterion writes
\begin{equation}\label{Mean-var-result}
\begin{split}
\Ur(\phi) &= R_0 + \left(\rho-c\right)\overline{\theta}T  +H[\phi-\Ir] \cdot \left(M(T)-cT\right)  - \gamma H\left[(\phi-\Ir)^2\right]\cdot M(T) \\
& - \gamma  H^2[\phi-\Ir] A(T) - \gamma H[\phi-\Ir] \cdot H[f\cdot  (\phi-\Ir)] \cdot B(T)   \\
\end{split}
\end{equation}
where
\begin{eqnarray*}
A(T) &=& 2 \int_0^T\int_0^te^{-\left(\beta- H[f]\right) (t-s)}\left[\beta\underline{\lambda}M(s) + m^{(2)}(s) \right]ds  \ dt- M^2(T) \ ,\\
B(T) &=& 2 \int_0^T\int_0^te^{-\left(\beta- H[f]\right)(t-s)}m(s)ds \ dt .
\end{eqnarray*}
\end{thm}
{\bf Proof of Theorem \ref{thmMeanVar}:}
The expectation of $R_T(\phi)$ can be deduced by Proposition \ref{Prop1}. We observe that in $R_T(\phi)$ given in \eqref{eqRTphi} there is only the third term as random, so obviously $\V\left[R_T(\phi)\right] =\V\left[X_T(\phi)\right]$. We have to compute the variance or, equivalently, the second order moment of $X_T[h]$, the computation of which is the result of the
following lemma. Given the Lemma, in order to obtain the result stated in the Theorem, is then sufficient to apply Lemma \ref{lemmaExQu} to $h= \phi-\Ir$. 
\hfill $\Box$

\medskip
\begin{lemma} \label{lemmaExQu}
The second order moment of $X_T[h]$ writes

\begin{equation}\label{eqEspX2}
\begin{split}
\E\left[X^2_T[h] \right] & =
H^2\left[h\right] \left(A(T) +M^2(T) \right)  \\
&+ H\left[h\right]\cdot H\left[f \cdot  h\right]\cdot B(T) \\
& + H\left[h^2\right] M(T) .
\end{split}
\end{equation}

\end{lemma}
{\bf Proof:} By It\^o formula (see \cite{Pro92} Theorem 32, p. 71), we have
\[X_T^2[h] = 2\int_0^TX_{t^-}[h] \int_0^{+\infty}h(z)\nu(dt,dz) + \int_0^T\int_0^{+\infty}h^2(z)\nu(dt,dz)\]
By taking the expectation, and using the definition of the compensator, we obtain
\begin{equation}\label{eqX2temp}
\E\left[X^2_T[h] \right] = 2 H[h] \int_0^T \E\left[ X_{t^-}[h]\; \lambda_{t^-} \right] dt + H\left[h^2\right]M(T)
\end{equation}
Now, let us focus on the term $\E\left[ X_{t^-}[h] \; \lambda_{t^-} \right]$. By It\^o formula, we have
\[\begin{aligned}
X_T[h]\cdot \lambda_T =&  \int_0^TX_{t^-}[h] \, d\lambda_t + \int_0^T\lambda_{t^-}dX_t[h] + \int_0^T\int_0^{+\infty}f(z) h(z) \nu(dt,dz)\\
= & \beta \int_0^TX_{t^-}[h] \left(\underline{\lambda}-\lambda_{t^-} \right)dt + \int_0^TX_{t^-} [h] \int_0^{+\infty} f(z)\nu(dt,dz)\\ 
& + \int_0^T\lambda_{t^-}\int_0^{+\infty}h(z)\nu(dt,dz) + \int_0^T\int_0^{+\infty}f(z)h(z)\nu(dt,dz) .
\end{aligned} \]
Let define for any $0 \le t \le T$ $$U(t) : = \E\left[X_t[h] \cdot \lambda_t\right].$$ By the definition of the compensator and equation \eqref{eqXh}, adapting the proof of Proposition \ref{Prop1} and using stochastic Fubini theorem \cite[Theorem IV.65, p 208]{Pro92}, we obtain the following (implicit) relation

\[
U(T) =H[h] \int_0^T \left(
\beta \underline{\lambda} M(t) 
+ m^{(2)}(t) \right) dt - \left(\beta-H[f]\right)\int_0^T U(t)dt + H[f\cdot h] \cdot  M(T) .
\]

We look at the solution under the form  $U(t)=e^{-\left(\beta- H[f]\right)t}\; V(t)$, with $V_0 = 0$ since by \eqref{eqXh} we have $U_0=0$. The function $V$ satisfies the following ordinary differential equation:
\[e^{-\left(\beta- H[f]\right) t}V^\prime(t) = H[h] \left( \beta\underline{\lambda} M(t) +  m^{(2)}(t) \right)+ H[f\cdot h]\cdot m(t)\]
Hence,
\begin{equation}\label{eqUT}
\begin{aligned}
U(T) = & H[h] \int_0^Te^{-\left(\beta-H[f]\right)(T-t)}\left(\beta\underline{\lambda} M(t) + m^{(2)}(t)\right) dt \\
 & + H[f\cdot  h ]\int_0^Te^{-\left(\beta-H[f]\right)(T-t)}m(t)dt 
\end{aligned}
\end{equation}

Using equations \eqref{eqX2temp}) and \eqref{eqUT}, we obtain the result. \hfill $\Box$

\medskip
The form of the mean-variance criterion is significantly simplified when the Hawkes process is reduced to a Poisson process as detailed in the following Corollary.

\begin{corollary}Assume that $\nu$ is Poisson process with (constant) intensity $\lambda_0$. Then the mean-variance criterion writes
\begin{equation}\label{criterium-Poisson}
\Ur(\phi)  = R_0 + (\rho-c) \overline{\theta}T + H[\phi-\Ir] \cdot \left( \lambda_0  - c  \right) T- \gamma H[(\phi-\Ir)^2]\lambda_0 T 
\end{equation}
That is only the first line in the equation \eqref{Mean-var-result}.
\end{corollary}

{\bf Proof:} The two extra terms in equation \eqref{Mean-var-result} vanishes in the Poisson setup. 
In fact Poisson case can be deduced from the Hawkes one by setting the self-exciting parameter $f$ and the mean-reverting speed $\beta$ to zero. The long run intensity $\underline{\lambda}$ will then coincide with $\lambda_0$. As a consequence we immediately have $M(t)= \lambda_0 t$, 
$A(t) = 0$, cancelling the first term into the second line of equation \ref{Mean-var-result}. The last term vanishes when the self exciting function $f$ cancels.
\hfill $\Box$

\medskip

We highlight that the expression \eqref{criterium-Poisson} of the mean-variance criterion in Poisson case depends on the loss only though the functions $(\phi-\Ir )$ and  $(\phi-\Ir)^2$.\\

\medskip
In the next corollary we show that, if the self exciting function is constant, i.e. $f(z) \equiv \overline{f}$, then the criterion also depends only on these two terms. As a consequence, the Hawkes case with constant feedback effect will produce the same optimal contract as in the Poisson case. Different optimal contracts can only arise when the feedback depends on the size of the loss. 

\begin{corollary}Assume that $f(z) \equiv \overline{f}$. Then, the mean-variance criterion writes
\begin{equation}\label{criterium-ConstantJump}
\begin{split}
\Ur(\phi) = & R_0 + \left(\rho-c\right)\overline{\theta}T +H[\phi-\Ir] \left(M(T)-cT\right)   \\
&- \gamma H^2[\phi-\Ir] \left(A(T)+ \overline{f} B(T)\right) - \gamma H\left[(\phi-\Ir)^2\right]M(T) 
\end{split}
\end{equation}
\end{corollary}

Formulae given by equations (\ref{criterium-Poisson}) and (\ref{criterium-ConstantJump}) do not depend on a specific term $H[f\cdot  (\phi-\Ir) ]$ as in the general case stated in Theorem \ref{thmMeanVar}. This additional term will give rise to specific form of reinsurance contracts as we will see in next section.

\section{Optimal Reinsurance Contracts}\label{secOptReinCont}

\subsection{Regularity of the criterion}

Let us introduce the following lemma on the regularity of the map $\Ur$. We show that  $\Ur$ admits a Fréchet-derivative; see for instance  Dieudonn\'e  \cite{Dieudonn} for classical definition of Fr\'echet-derivative of a map defined on a normed space.

\begin{lemma} The mapping $\Ur : L^2(\Theta) \longrightarrow \R$  is Fr\'echet-derivable in $L^2(\Theta)$ and the  Fr\'echet-derivative is a function in $L^2(\Theta)$ 
given by $D\Ur(\phi) (z)= G(\phi)\Theta(dz)$ with
\[\begin{split}
G(\phi) := & \Big( (M(T) - cT)  - 2\gamma H[\phi-\Ir] A(T) - \gamma H[f\cdot (\phi-\Ir)]B(T)\Big) \mathds{1}  \\
& -\gamma  H[\phi-\Ir]B(T)\;  f -2\gamma M(T)\; \left( \phi-\Ir\right)
\end{split}\]
where $\mathds{1} $ denotes the constant function.
\end{lemma}

{\bf Proof}. 
We start from the expression of $\Ur(\phi)$ given by Theorem \ref{thmMeanVar}.
The equation \eqref{Mean-var-result} shows that $\Ur$ depends on $\phi$ through three terms: $H[\phi-\Ir]$, $H[(\phi-\Ir)^2]$ and $H[f \cdot (\phi-\Ir)]$. We will study the  Fr\'echet-derivative of the three terms. The final result will then be obtained as the sum of three Fr\'echet-derivative terms. 
Firstly, map $H[\phi-\Ir]$ is continuous and linear on $L^2(\Theta)$, hence Fr\'echet-derivable, with  Fréchet-derivative constant given by $g \mapsto \int_{0}^{+\infty}g(z)\Theta(dz)$. 
 Now, let us turn to the term $H[(\phi-\Ir)^2]$. We have, for any $g \in L^2(\Theta)$,
\[\begin{aligned}
H[(\phi- \Ir +g)^2] =&   \int_0^{+\infty} \left(\phi(z)-z\right)^2\Theta(dz) + 2 \int_0^{+\infty}\left(\phi(z)-z\right)g(z)\Theta(dz) \\ &+ \int_0^{+\infty}g^2(z)\Theta(dz).
\end{aligned}\]
Then, the Fréchet-derivative of $H[(\phi - \Ir)^2]$ is $g \mapsto 2\int_{0}^{+\infty}g(z)\left(\phi-\Ir\right)\Theta(dz)$. 
Finally, we focus on $H[f \cdot (\phi-\Ir)]$. We have, for any $g \in L^2(\Theta)$,
 $$
 H[f\cdot (\phi- \Ir+g) ] = \int_0^{+\infty}f(z)\left(\phi(z)-z\right)\Theta(dz) + \int_0^{+\infty}f(z)g(z)\theta(dz).$$
Then, the Fr\'echet-derivative of $H[f \cdot (\phi-\Ir)]$ is $g \mapsto \int_{0}^{+\infty}g(z)f(z) \Theta(dz)$.

 \hfill $\Box$

\bigskip
We notice that, in the constant case i.e. $f\equiv \overline{f}\mathds{1}$, the Fréchet-derivative depends only on $H[\phi-\Ir]$, which is what we have in the usual Poisson case. \\

So, from now on, we assume the linear case i.e. $f(z) = \Lambda z$.

\subsection{The Solution of the Program} 

\begin{proposition} The insurer program given by
\[\left(\Pr\right) : \:\: \sup_{\phi \in \Cr}\Ur(\phi)\]
admits a solution $\phi^{\ast} \in \Cr$
\end{proposition}
{\bf Proof}. The set $\Cr$ is closed, convex and bounded in $L^2(\Theta)$. It is therefore a compact for the weak topology (cf. Corollary 3.22, p. 71 in \cite{Bre}). The function $\Ur$ is derivable, hence continuous for the norm topology on $L^2(\Theta)$, thanks to the previous lemma. Hence, it admits a minimum on $\Cr$. \hfill $\Box$

\medskip

We can now state the main result of our paper:

\begin{thm}\label{thmSolution} Assume the support of $\Theta(dz)$ be not bounded, $\Lambda>0$ and $c T - M(T) >0$. Then the optimal solution $\phi^{\ast}$ of the insurer program verifies the following.
\begin{enumerate}
\item On $\Ar_0 := \left\{z \vert G(\phi^{\ast})(z) <  0 \right\}\not=\emptyset$, $\phi^{\ast}\equiv 0$.
\item On $\Ar_1 := \left\{z \vert G(\phi^{\ast})(z) =0 \right\}\not=\emptyset$, the optimal contract is affine with a slope larger than $1$ such that the the following implicit system is satisfied
\begin{equation}\label{eqA1}
\phi^{\ast}(z) = \left( 1 - \Lambda\frac{B(T)}{2 M(T)} \, H[\phi^{\ast}-\Ir] \right) z - C^{\ast} 
\end{equation}
with
\[C^{\ast}:= \frac{cT-M(T) +2\gamma A(T)\, H[\phi^{\ast}-\Ir] +\gamma B(T)  H[\Lambda \Ir \cdot (\phi^{\ast}-\Ir)] }{2\gamma M(T)}.\] 

\item On $\Ar_2 := \left\{z \vert G(\phi^{\ast})(z) > 0 \right\}\not=\emptyset$, $\phi^{\ast}(z) = z$.
\end{enumerate}
\end{thm}
\begin{rem}
Assumption $c T - M(T) >0$  means that the reinsurance cost is larger than the expected cumulated intensity over $[0,T]$.
\end{rem}
{\bf Proof}. We follow the approach of Raviv \cite{Raviv} based on variational inequalities. An approach by penalisation has also been proposed by Gollier \cite{Gollier}.

Let $\phi \in \Cr$ and $\alpha \in [0,1]$. We have $\Ur\left(\phi^{\ast}\right) \ge \Ur\left(\phi^{\ast} + \alpha(\phi -\phi^{\ast}) \right)$. Letting $\alpha$ go to $0$, we obtain the first order necessary condition
\begin{equation}\label{eqCN1}
\forall \phi \in \Cr, \:\: \int_0^{+\infty}G(\phi^{\ast})[z]\left(\phi(z) - \phi^{\ast}(z) \right) \le 0
\end{equation}
First of all, assuming by contradiction that the optimal solution is full reinsurance, i.e. $\phi^{\ast}=\Ir$ , we would have $H[\phi^{\ast}-\Ir]=H[\Lambda \Ir \cdot (\phi^{\ast}-\Ir)]=0$ and  $G(\phi^{\ast})  = (M(T) - cT)<0$. It is in contradiction with equation \eqref{eqCN1}. Hence, the optimal solution is different from the total reinsurance: $\phi^{\ast}\not=\Ir$. 

Second, assuming by contradiction that the optimal solution is no-reinsurance at all, i.e. $\phi^{\ast}=0$, we have $H[-\Ir]= -\overline\theta <0$.
A direct computation gives
 $ G(0)= (M(T) - cT)  + 2\gamma \overline\theta A(T) + \gamma H[\Lambda \Ir^2]B(T) + \gamma \left[\overline{\theta}\Lambda B(T) + 2M(T)\right] \Ir $, for $\Theta$-almost surely every $z$. This result violates condition \eqref{eqCN1}, since the support of $\Theta$ is not bounded.

Let $\varepsilon$ small enough and $D_{\varepsilon}:=\left\{z \vert \varepsilon < \phi^{\ast} < z-\varepsilon \right\}$, which is non empty thanks to the assumption and the previous step of the proof. Then, $\phi^{\ast}_{+}:= \phi^{\ast}+\varepsilon \indic{D_{\varepsilon}}$ and $\phi^{\ast}_{-}:= \phi^{\ast}-\varepsilon \indic{D_{\varepsilon}}$ are in $\Cr$. It yields $G(\phi^{\ast})=0$ on   $D_{\varepsilon}$, thanks to \eqref{eqCN1}.
Taking the limit of $\varepsilon$ goes to zero, we obtain  $G(\phi^{\ast})=0$ on $\Ar_1$  $\Theta$-almost surely thanks to monotone class theorem. Direct computation of this condition yields $\Ar_1$.
The two other domains, $\Ar_0$ and $\Ar_1$ are given by the same kind of arguments.

Now, let us characterise the sets $\Ar_0$, $\Ar_1$ and $\Ar_2$. First, we notice that the gradient of $\Ur$ is increasing on $\Ar_0$ and the slope of $\phi^{\ast}$, as given in equation \eqref{eqA1}, is greater than $1$ (because $H[\phi^{\ast} -\Ir]<0$). Therefore, we have $\Ar_0 = [0,a]$, $\Ar_1 = ]a,b]$ and $\Ar_2 = ]b,+\infty[$. We also deduce from that configuration that $a$ cannot be equal to $0$. Indeed, assume that $a=0$, then $\phi^{\ast}(z) > z$ because its slope is larger than $1$ on $\Ar_1$, which is not possible. In this case, the only solution would be $\Ar_1 = \emptyset$, which is not the case from the previous steps of the demonstration. Now, we can see that $\Ar_2$ cannot be empty because the slope of $\phi^{\ast}$ is larger than $1$ on $\Ar_1$. \hfill $\Box$

\medskip

\begin{rem}
The form of the solution given by Theorem \ref{thmSolution} is displayed in Figure \ref{fig:OptimalSolution}. We can see that $\phi^{\ast}$ is piecewise affine, with three different parts. When the solution is not trivial (in the interval $[a,b]$), its slope is greater than $1$ and the extra-slope is driven by $\Lambda$, that is
the self-exciting parameter. This part stems from the clusters effects through the term $B(T)$. This slope greater than one induces more reinsurance than in the classical case where the deductible is optimal. This leads to a zone where the extreme risks are totally reinsured. Corollary \ref{corFormSolution} specifies the form of the optimal solution as a function of $a$ and $b$. 
\end{rem}
                                                                                                                         
 \begin{corollary}\label{corFormSolution} Let the assumptions of Theorem \ref{thmSolution} hold. Then, the optimal solution of $\left(\Pr\right)$ is given by
\[\phi^{\ast}(z) = \min\left\{z\,; \,\frac{b}{b-a}(z-a)_+ \right\}\]
with $b>a$ satisfying
\[\Lambda \int_0^b\left(z - \frac{b}{b-a}(z-a)_+\right)\Theta(dz) = \frac{2  M(T)}{B(T)}\frac{a}{b-a} .\]

\end{corollary}
{\bf Proof}. From Theorem \ref{thmSolution}, we know that $\phi^{\ast}$ is affine on $[a,b]$, equal to $0$ at $a$ and equal to $b$ at $b$. Its slope is, therefore, $\frac{b}{b-a}$. It is also equal, by equation \eqref{eqA1}, to $1 - \Lambda\frac{B(T)}{2 M(T)} H[\phi^{\ast} -\Ir]$. Direct computation yields the result. \hfill $\Box$

\medskip
\subsection{Limit to the Poisson case}

In this subsection, we show how the previous result converges to the Poisson framework when the 
self-exciting parameter $\Lambda$ goes to zero. 
Since the parameter $\Lambda$ impacts not only the shape of the optimal contract but also the frequency of the losses over the considered window time $[0,T]$, we will assume that the parameters are changed in such way that the following holds.

\begin{description}
\item[Invariance of the expected number of events:] the frequency of the loss arrival in the Poisson case, denoted by $\lambda_P$ coincide with the average intensity of the Hawkes arrival. That is $\lambda_P = M(T)/T$.
\item[Invariance of the expected cost:] The optimal contract is specified by its cost $c H[\phi^\ast]$. We will assume that the cost is constant and then   $H[\phi^\ast]$ is invariant when $\Lambda$ changes.
\end{description}

The first condition guarantee that the average number of events during the window $[0,T]$ is unchanged 
when $\Lambda$ decreases. A direct computation, adapting for instance the arguments in \cite[Section 5.1]{BerSalSco18},  the  average number of events $M(T)$ increases with the self-exciting parameter $\Lambda$. We have, in particular, $M(T) = \frac{\beta}{\beta - \overline{\theta} \Lambda} \underline{\lambda} T + (\lambda_0 - \underline{\lambda}) e^{- (\beta - \overline{\theta} \Lambda) T}$.
As a consequence, when $\Lambda$ decreases, we have to increase the parameter $\underline{\lambda}$ to keep $M(T)= \lambda_P T$. 

To simplify our analysis, we assume that the contract is signed outside cluster period, in order to can assume $\lambda_0 = \underline{\lambda}$.
We sum-up the results into the following corollary.

\begin{corollary}[Convergence to Poisson setup]
Under the previous conditions, the optimal contract $\phi^\ast$ converges to the excess-of-loss reinsurance contract when the self-exciting parameter $\Lambda$ decreases to zero and in particular we have that 
\begin{description}
\item[Decreasing slope:] The slope in the intermediate region $\Ar_1$ decreases when $\Lambda$ decreases and converges to $1$ if $\Lambda$ goes to zero.
\item[Decreasing franchise:] The deductible threshold $a$ decreases when $\Lambda$ decreases, such that the optimal contract is to not cover any loss 
below to $a$.
\item[Increasing threshold $b$:] The threshold $b$ increases when $\Lambda$ decreases, such that the optimal contract is to cover all loss 
above $b$.

\end{description}
\end{corollary}

{\bf Proof}.
Since $\lambda_0 = \underline{\lambda}$, we directly obtain $m(t) = \frac{\beta}{\beta - \overline{\theta} \Lambda} \underline{\lambda} = \lambda_P$, see \cite[Section 5.1]{BerSalSco18}, and then $m(t)$ is
unchanged when $\Lambda$ changes. By a direct integration we have
$$
B(T) = \frac{1 - e^{-(\beta - \Lambda \overline{\theta}) T}}{(\beta - \Lambda \overline{\theta})^2 }\lambda_P 
$$
It is easy to see that $B(T)$ is increasing in $\Lambda$. It is easy then to see that $\Lambda B(T)$ is a strictly increasing function of $\Lambda$. When $\Lambda$ decreases the slope of the optimal contract 
in the intermediate region $\Ar_1$ is then decreasing.

The two other results are obtained by contradiction. Assume first that there exists $\Lambda_1>\Lambda_2$ such that $a(\Lambda_1) < a(\Lambda_2)$. Then, since the slope of optimal contract in the intermediate region $\Ar_1$ increases with $\Lambda$, it is easy to see that the optimal contract for $\Lambda_1$ dominates the one for  $\Lambda_2$ violating the condition  of invariant cost $c H[\phi^\ast]$.
A similar argument guarantees that the threshold $b$ increases when $\Lambda$ decreases.
\hfill $\Box$

\medskip

Finally, we can easy remark that the mean variance criterion of the Hawkes setup, see equation \eqref{Mean-var-result}, converges to the one of Poisson, see \eqref{criterium-Poisson}, since the last term is proportional to $\Lambda$ whereas the second-to-last is proportional to $A(T)$ which disappears when $\Lambda$ goes to zero. In particular the function $A(T)$ can be interpreted as a variance of the number of the jumps expected in the window $[0,T]$. In this sense the second-to-last term in \eqref{Mean-var-result} can be interpreted as the result of a Wald identity that will disappear in the Poisson limit.

\section{Conclusion}

In a static framework, the optimal form of reinsurance contract is the excess-loss contract, as demonstrated by Arrow \cite{Arr}. This result has been extended to a continuous-time framework \cite{Aas}. However, we show that this is no longer valid in the presence of clustering effects on losses. We consider a marked Hawkes process, where the times of jumps represent the occurrence of claims, and the marks indicate the associated financial losses. The intensity of the Hawkes process increases with each jump, reflecting the presence of clusters of events. The reinsurance contract is modelled as a function in $L^2(\R^+)$
and applies to each loss occurrence. The criterion used to assess the risk profile is a static mean-variance function over a given finite horizon $T$.
We derive a closed-form formula for this criterion and demonstrate that the optimal contract is not an excess-loss function (deductible), but rather comprises three components:
\begin{itemize}
\item No reinsurance below a certain level of loss $a$
\item Reinsurance according to an affine function on $[a, b]$, with a slope greater than $1$. This slope, which exceeds $1$, is a consequence of clustering effects and stochastic jump intensity, and it accounts for the existence of the third component.
\item Full reinsurance above the level of loss $b$
\end{itemize}
This result is derived under very mild assumptions regarding the distribution of loss sizes (with non-bounded support) and the reinsurance cost (which exceeds the expected intensity of jumps). The specific form of the optimal contract in the clustered case is a direct consequence of the stochastic jump intensity, which influences the variance. When our model converges to a Poisson process while maintaining a constant average intensity over the interval  $[0,T]$, the optimal contract aligns with the classical excess-loss problem.

This result highlights that clustering effects can significantly alter the perception of risk and the need for reinsurance. Specifically, the optimal form of reinsurance in the presence of clusters of jumps tends to fully cover tail risk.
We emphasize that, in contrast to a substantial body of literature \cite{HarRav, HubMaySmi}, our setup does not lead to a moral hazard issue. First, ex-ante moral hazard, referring to the distortion of insurers' behaviour, does not arise due to the feedback effect of the Hawkes process, which increases the probability of new events. Meanwhile, ex-post moral hazard, related to the distortion of claims, is canceled in our marked Hawkes setup, where the arrival of new events is positively influenced not only by previous occurrences but also by their magnitudes.

Future investigations could include a statistical study of the feedback self-exciting effect, in particular relaxing the hypothesis of exponential kernel using Laplace transform technique see \cite{BacMuz, BLP24, BPS22, BouCherHil1, HilRevRos}. Techniques proposed by Brignone et al. \cite{BriGonSga} in the context of commodities could be adapted for this purpose. Additionally, Brachetta et al. \cite{BraCalCecSga} have examined a mixture of Poisson and Hawkes frameworks without focusing on the shape of the optimal contract, providing an opportunity to analyse how the optimal contract depends on the chosen filter. Finally, a more general form of contracts could be explored, such as relaxing the non-negative constraint as suggested by Gollier \cite{Gollier}.

\section*{Acknowledgement}
Cristina Di Girolami acknowledges the support of INDAM and GNAMPA, in particular Project CUP E53C23001670001 \emph{Problemi di controllo ottimo stocastico con memoria a informazione parziale}. Cristina Di Girolami acknowledges the support of Ministero Italiano dell'Università e della Ricerca, through the Programma di Ricerca Scientifica
di Rilevante Interesse Nazionale 2022BEMMLZ\_003  \emph{Stochastic control and games and the role of information}. 
Simone Scotti acknowledges the support of University of Pisa, under the PRA - Progetti di Ricerca di Ateneo and the support of Institut Louis Bachelier - Institute Europlace of Finance.


\begin{thebibliography}{99}

\bibitem{Aas} Aase, K. K. (1993): Premiums in a dynamic model of a reinsurance market. Scandinavian Actuarial Journal, 1993(2), 134-160. 

\bibitem{AberJedi}
Abergel, F.,  Jedidi, A. (2015). Long-time behavior of a Hawkes process--based limit order book. SIAM Journal on Financial Mathematics, 6(1), 1026-1043.

\bibitem{AbiLarPul}
Abi Jaber, E., Larsson, M.,  Pulido, S. (2019). Affine Volterra processes. The Annals of Applied Probability, 29(5), 3155-3200.

\bibitem{AlAs06}
Albrecher, H., Asmussen c, S. R. (2006). Ruin probabilities and aggregate claims distributions for shot noise Cox processes. Scandinavian Actuarial Journal, 2006(2), 86-110.

\bibitem{Arrow} Arrow, K. J. (1963). Uncertainty and the Welfare Economics of Medical Care. The American Economic Review, 53(5), 941-973.

\bibitem{Arr} Arrow, K.J. (1974): Optimal insurance and Generalized Deductibles. Scandinavian Actuarial Journal 1, 1-42.

\bibitem{BacMuz}
Bacry, E.,  Muzy, J. F. (2014). Hawkes model for price and trades high-frequency dynamics. Quantitative Finance, 14(7), 1147-1166.

\bibitem{BalGhe}
Baldwin, A., Gheyas, I., Ioannidis, C., Pym, D., \& Williams, J. (2017). Contagion in cyber security attacks. Journal of the Operational Research Society, 68(7), 780-791.

\bibitem{Barsotti}
Barsotti, F., Milhaud, X., \& Salhi, Y. (2016). Lapse risk in life insurance: Correlation and contagion effects among policyholders? behaviors. Insurance: Mathematics and Economics, 71, 317-331.

\bibitem{Bauerle} Bauerle, N. (2005). Benchmark and mean-variance problems for insurers. Mathematical Methods of Operations Research, 62, 159-165.

\bibitem{Ber02} Bernis, G. (2002): Equilibrium in a reinsurance market with short sale constraints. { Economic Theory}, 20, p. 295-320.

\bibitem{BGSS23}
Bernis, G., Garcin, M., Scotti, S.,  Sgarra, C. (2023). Interest rates term structure models driven by Hawkes processes. SIAM Journal on Financial Mathematics, 14(4), 1062-1079.

\bibitem{BerSalSco18} Bernis, G., Salhi, K., and Scotti, S. (2018): Sensitivity analysis for marked Hawkes processes: application to CLO pricing. Mathematics and Financial Economics, 12(4), 541-559.

\bibitem{BerSco20} Bernis, G.  and Scotti, S. (2020): Clustering Effects via Hawkes Processes, in {\it From Probability to Finance}, Jiao Editor, Springer.

\bibitem{BLP24}
Bondi, A., Livieri, G., Pulido, S. (2024). Affine Volterra processes with jumps. Stochastic Processes and their Applications, 168, 104264.

\bibitem{BPS22}
Bondi, A., Pulido, S.,  Scotti, S. (2024). The rough Hawkes Heston stochastic volatility model. Mathematical Finance, 34, 1197-1241.


\bibitem{Bor} Borch, K. H. (1962): Equilibrium in a reinsurance market. Econometrica 30, 424–444

\bibitem{BouCherHil1} Bessy-Roland Y., Boumezoued, A. and  Hillairet, C. (2021) Multivariate Hawkes process for cyber insurance. 
Annals of Actuarial Science , Volume 15 , Issue 1 , March 2021 , pp. 14 - 39. 

\bibitem{BraCec}
Brachetta, M., Ceci, C. (2019): Optimal proportional reinsurance and investment for stochastic factor models.
Insur. Math. Econ. 87, 15-33 

\bibitem{BraCec2}
 Brachetta, M., Ceci, C. (2020): A BSDE-based approach for the optimal reinsurance problem under partial
information. Insur. Math. Econ. 95, 1-16 

\bibitem{BraCalCecSga} Brachetta, M., Callegaro, G., Ceci, C. and Sgarra, C. (2024) Optimal reinsurance via BSDEs in a partially observable model with jump clusters. Finance and Stochastics 28, 453-495.

\bibitem{Bre} Brezis, H. (2010): Functional Analysis, Sobolev Spaces and Partial Differential Equations. Springer, Heidelberg.

\bibitem{BriGonSga}
Brignone, R., Gonzato, L.,  Sgarra, C. (2024). Hawkes Processes in Energy Markets: Modelling, Estimation and Derivatives Pricing. In Quantitative Energy Finance: Recent Trends and Developments (pp. 41-72). Cham: Springer Nature Switzerland.

\bibitem{BriSga20} Brignone, R., and Sgarra, C. (2020): Asian options pricing in Hawkes-type jump-diffusion models. Annals of Finance, 16(1), 101-119.

\bibitem{CalMazSga}
Callegaro, G., Mazzoran, A.,  Sgarra, C. (2022). A self-exciting modeling framework for forward prices in power markets. Applied stochastic models in business and industry, 38(1), 27-48.


\bibitem{CallFonHillOng}
Callegaro, G., Fontana, C., Hillairet, C., Ongarato, B. (2025) A Stochastic Gordon-Loeb model for cybersecurity investment under clustered attacks. Submitted.



 \bibitem{DassJan}  Dassios, A., Jang, J.W. (2003): Pricing a catastrophe reinsurance and derivatives using the Cox process with
shot noise intensity. Finance Stoch. 7, 73-95 

 \bibitem{DassJan2}  Dassios, A., Jang, J.W. (2005): Kalman-Bucy filtering for linear systems driven by the Cox process with shot
noise intensity and its application to the pricing of reinsurance contracts. J. Appl. Probab. 42, 93-107

 \bibitem{DassZh}
Dassios, A., Zhao, H. (2011): A dynamic contagion process. Adv. Appl. Probab. 43, 814-846

 \bibitem{Dieudonn}
Dieudonn\'e, Jean (1969), Foundations of modern analysis, Boston, MA.

\bibitem{EecGolSch} Eeckhoudt, L., Gollier, C and Schlesinger, H.  (1996): Changes in Background Risk and Risk Taking Behavior. Econometrica, 64(3), 683-689.

\bibitem{EmSchGra}
Embrechts, P., Schmidli, H., Grandell, J.(1993): Finite-time Lundberg inequalities in the Cox case. Scand.
Actuar. J. 1, 17-41 

\bibitem{GarLeoTor}
Garetto, M., Leonardi, E., Torrisi, G. L. (2021). A time-modulated Hawkes process to model the spread of COVID-19 and the impact of countermeasures. Annual reviews in control, 51, 551-563.

\bibitem{Gollier}
Gollier, C. (1987). The design of optimal insurance contracts without the nonnegativity constraint on claims. The Journal of Risk and Insurance, 54(2), 314-324.

\bibitem{GolSch} Gollier, C., Schlesinger, H.  (1996): Arrow's Theorem on Optimality of Deductibles: A Stochastic Dominance Approach. Economic Theory 7, 359-363

\bibitem{Grandell}
Grandell, J. (1991): Aspects of Risk Theory. Springer, New York.

\bibitem{HarRav}
Harris, M., Raviv, A. (1978). Some results on incentive contracts with applications to education and employment, health insurance, and law enforcement. The American economic review, 68(1), 20-30.


\bibitem{Hawkes}
Hawkes, A.G. (1971): Spectra of some self-exciting and mutually exciting point processes. Biometrika 58,
83-90.

\bibitem{HawOak}
Hawkes, A. G.,  Oakes, D. (1974). A cluster process representation of a self-exciting process. Journal of applied probability, 11(3), 493-503.



\bibitem{HilRevRos}  Hillairet, C.,  Reveillac, A.,   Rosenbaum, M. (2023) An expansion formula for Hawkes processes and application to cyber-insurance derivatives, Stochastic Processes and their Applications, Volume 160,  89-119.

\bibitem{HorstXu}
Horst, U.,  Xu, W. (2019). A scaling limit for limit order books driven by Hawkes processes. SIAM Journal on Financial Mathematics, 10(2), 350-393.

\bibitem{HorstXu2}
Horst, U.,  Xu, W. (2022). The microstructure of stochastic volatility models with self-exciting jump dynamics. The Annals of Applied Probability, 32(6), 4568-4610.

\bibitem{HubMaySmi}
Huberman, G., Mayers, D., Smith Jr, C. W. (1983). Optimal insurance policy indemnity schedules. The Bell Journal of Economics, 415-426.

\bibitem{JiZh}
Jiang, H.,  Zhang, Z. (2025). Equity-linked annuity valuation under fractional jump-diffusion financial and mortality models. Scandinavian Actuarial Journal, 300-339.

\bibitem{JMS}
Jiao, Y., Ma, C., Scotti, S. (2017). Alpha-CIR model with branching processes in sovereign interest rate modeling. Finance and Stochastics, 21(3), 789-813.

\bibitem{JMSZ}
Jiao, Y., Ma, C., Scotti, S., Zhou, C. (2021). The Alpha-Heston stochastic volatility model. Mathematical finance, 31(3), 943-978.

\bibitem{KarShr91} Karatzas, I. and Shreve, S. (1991): Brownian Motion and Stochastic Calculus ($2^{nd}$ edition). Springer, Heidelberg.

\bibitem{LesDeaLej} Lesage, L., Deaconu, M., Lejay, A. et al. (2022) Hawkes Processes Framework With a Gamma Density As Excitation Function: Application to Natural Disasters for Insurance. Methodol Comput Appl Probab 24, 2509-2537.

\bibitem{LiaBay}
Liang, Z., Bayraktar, E. (2014): Optimal reinsurance and investment with unobservable claim size and
intensity. Insur. Math. Econ. 55, 156-166


\bibitem{Pro92} Protter, Ph. (1992): Stochastic Integration and Differential Equations. Springer, Heidelberg.

\bibitem{Raviv}
Raviv, A. (1979). The Design of an Optimal Insurance Policy. The American Economic Review, 69(1), 84-96.


\bibitem{Touzi} Touzi, N. (2000). Optimal insurance demand under marked point processes shocks. Annals of Applied Probability, 283-312.

\bibitem{Yan}
Yan, J. A. (1982). A propos de l'int\'egrabilit\'e uniforme des martingales exponentielles. S\'eminaire de probabilit\'es de Strasbourg, 16, 338-347.



\end{thebibliography}
\end{document}